\documentclass[a4paper,10pt,reqno]{amsart}

\usepackage[utf8]{inputenc}
\usepackage[english]{babel}

\usepackage[margin=40mm,top=39mm,bottom=39mm]{geometry}
\usepackage{amsmath,amsthm,amsfonts,amssymb}

\usepackage{hyperref}
\hypersetup{
  colorlinks = true,
pdfstartview = FitH,
 pdfpagemode = UseNone,
   pdfauthor = {Carofiglio Cherubini Gambini},
  pdfcreator = {Pdflatex},
}
\numberwithin{equation}{section}
\numberwithin{figure}{section}
\numberwithin{table}{section}
\makeatletter
\let\c@table\c@figure
\let\ftype@table\ftype@figure
\makeatother
\usepackage{tikz}
\usepackage{colortbl}
\usepackage[size=small]{caption}
\usepackage{cite}

\theoremstyle{plain}
\newtheorem{theorem}{Theorem}[section]
\newtheorem{obs}[theorem]{Observation}

\theoremstyle{definition}
\newtheorem{conj}[theorem]{Conjecture}

\author[Carofiglio]{Leonardo Carofiglio}
\address{
         Dipartimento di Matematica Guido Castelnuovo
         Sapienza Universit\`a di Roma
         Piazzale Aldo Moro, 5
         00185 Rome, Italy
        }
\email{carofiglio.1916715@studenti.uniroma1.it}

\author[Cherubini]{Giacomo Cherubini}
\address{
        Istituto Nazionale di Alta Matematica ``Francesco Severi'',
        Research Unit Dipartimento di Matematica ``Guido Castelnuovo'',
        Sapienza Universit\`a di Roma, Piazzale Aldo Moro 5, I-00185, Roma
        }
\email{cherubini@altamatematica.it}

\author[Gambini]{Alessandro Gambini}
\address{
         Dipartimento di Matematica Guido Castelnuovo
         Sapienza Universit\`a di Roma
         Piazzale Aldo Moro, 5
         00185 Rome, Italy
        }
\email{alessandro.gambini@uniroma1.it}

\title[Title]{On Eswarathasan--Levine and Boyd's conjectures for harmonic numbers}

\keywords{harmonic number, harmonic sum, Wolstenholme prime, harmonic prime}
\subjclass[2020]{Primary 11B83; Secondary 11Y55, 11Y70}

\begin{document}

\begin{abstract}
We provide numerical evidence towards three conjectures on harmonic numbers
by Eswarathasan--Levine and Boyd. Let $J_p$ denote the set of integers 
$n\geq 1$ such that the harmonic number $H_n$ is divisible by a prime $p$.
The conjectures state that:
$(i)$ $J_p$ is always finite and of the order $O(p^2(\log\log p)^{2+\epsilon})$;
$(ii)$ the set of primes for which $J_p$ is minimal
(called harmonic primes) has density $e^{-1}$ among all primes;
$(iii)$ no harmonic number is divisible by $p^4$.
We prove $(i)$ and $(iii)$ for all $p\leq 16843$ with at most one exception,
and enumerate harmonic primes up to~$50\cdot 10^5$, finding a proportion
close to the expected density. Our work extends previous computations
by Boyd by a factor of about $30$ and $50$, respectively.
\end{abstract}

\maketitle\thispagestyle{empty}

\section{Introduction}

The sequence of harmonic numbers
\[
H_n=1+\frac12+\frac13+\cdots+\frac1n
\]
is a much studied one in the literature, mainly due
to its connections with the Riemann zeta function and Bernoulli numbers.
Over the centuries, many arithmetic properties of $H_n$
have been discovered; a well-known example is Wolstenholme’s theorem \cite{Wolstenholme},
which states that $p^2$ divides the numerator of $H_{p-1}$
for every prime $p\geq 5$. More generally, the divisibility
of $H_n$ by a given prime $p$ has always attracted a lot of interest
\cite{A,Bayat,BHMT,Boyd,Carlitz,EL,Gl1,Gl2,HT,Hong,S,SunHong,WC}.

Our motivation for looking into the sequence $H_n$ is threefold.
First, harmonic numbers are related to $p$-adic $L$-functions \cite{W}, which are less
well-understood than the classical ones.
A striking fact in this context is that we do not even know
if $p$-adic zeta functions are always non-zero on the positive integers
(see e.g~\cite{Beukers,Cal}).

Second, the set of harmonic numbers divisible by a given prime $p$
can be described by a probabilistic model, which allows one to make
conjectures on what we should expect. This has been worked out in full
detail by Boyd~\cite{Boyd}.

Third, specialized software is available to test the predictions
made by the probabilistic model. As explained by Boyd in \cite[\S 5]{Boyd},
the `naive' approach of computing $H_n$ from $H_{n-1}$
and then checking the divisibility is unfeasible for large values of $n$.
Instead, a better-tailored $p$-adic method can reach much higher values.

\subsection{The set \texorpdfstring{$J_p$}{J\_p}}
The central object in this paper is the set
\[
J_p := \{n\geq 1:\; \nu_p(H_n)\geq 1\},
\]
where $p$ is a prime and $\nu_p(a)$ denotes the $p$-adic valuation of $a$.
In other words, $J_p$ contains those $n$ such that $p$
divides the numerator of $H_n$ (when written in lowest terms).
We aim to describe the two extreme cases of how small and how large
the cardinality $|J_p|$ can be.

In 1991, Eswarathasan and Levine~\cite{EL} initiated
a study of $J_p$ and computed the full set when $p=3,5,7$.
Based on the fact that these sets were all finite,
they conjectured that this should always be the case
(a rather ambitious conjecture, viewed the limited evidence).

\begin{conj}\label{conj1}
The set $J_p$ is finite for all primes $p$.
\end{conj}

In the same paper, they showed that for all $p\geq 5$
the set $J_p$ always contains $p-1$, $p^2-p$ and $p^2-1$.
They called \emph{harmonic primes} those $p$ for which $|J_p|=3$
and suggested that they should occur infinitely often.

\begin{conj}\label{conj2}
There are infinitely many harmonic primes.
\end{conj}

To explore these conjectures, Eswarathasan and Levine devised an algorithm based on the decomposition
\begin{equation}\label{1502:eq001}
H_{pn+k} = H_{pn+k}^* + \frac{H_{n}}{p},
\end{equation}
where $k\in[0,p-1]$ and $H_n^*$ denotes a sum as in $H_n$ but restricted to integers coprime to $p$.
Since $H_{pn+k}^*\equiv H_{k}$ modulo $p$, it follows from \eqref{1502:eq001} that
\begin{equation}\label{1502:eq002}
H_{pn+k} \equiv H_{k} + \frac{H_n}{p} \pmod{p}.
\end{equation}
Therefore, $pn+k\in J_p$ if and only if we have
$n\in J_p$ and $p^{-1}H_n\equiv -H_k$ modulo~$p$
(cf. \cite[Theorem~3.1]{EL}).
In particular, this suggests a search strategy as follows:
after computing $H_k$ modulo~$p$ for all $k=1,\dots,p-1$,
determine the elements of $J_p\cap [p^m,p^{m+1}-1)$ and then use the above criterion
to find $J_p\cap[p^{m+1},p^{m+2}-1)$.

In 1994, Boyd~\cite{Boyd} extended this method by exploiting
a $p$-adically convergent series for $H_{pn}-p^{-1}H_n$ (see~\cite[Theorem~5.2]{Boyd}),
which allowed him to essentially iterate the recursion in \eqref{1502:eq002}
and get back to computing only the initial interval $J_p\cap[1,p-1]$, but to a high
$p$-adic precision. He managed to establish that $J_p$ is finite
for all primes $p \le 547$ except possibly for $p\in\{83, 127, 397\}$.

Boyd also explained how the set $J_p$ can be described
in terms of a probabilistic Galton--Watson branching process.
Such a random model suggests that, with probability one, the cardinality
$|J_p|$ is indeed finite (in agreement with Conjecture \ref{conj1})
and of the order $O(p^2(\log\log p)^{2+\epsilon})$,
with infinitely many primes satisfying $|J_p|\geq p^2(\log\log p)^2$.
In addition, Boyd's model predicts that harmonic primes should have
density $e^{-1}$ among all primes, which gives a quantitative refinement
of Conjecture~\ref{conj2}. Finally, it predicts that $H_n$
cannot be divisible by high powers of $p$ \cite[p.288]{Boyd}.

\begin{conj}\label{conj3}
There are no pairs $(p,n)$ with $\nu_p(H_n)\geq 5$.
The case $\nu_p(H_n)=4$, if it occurs at all, should occur finitely many times.
\end{conj}

In contrast, it is very common that $\nu_p(H_n)\leq 2$.
The case $\nu_p(H_n)=3$ occurs, too, although rarely,
the first instance being when $p=11$ and $n=848$.

\subsection{Main result}
We extend Boyd's results in two directions.
First, in a `vertical direction', so to speak, we consider small primes and check how
large $|J_p|$ can get. For a single prime, this can become very time-consuming
and so we decided to stop at $p=16843$ (the first Wolstenholme prime),
extending Boyd's computations by a factor of about~30.
In a `horizontal direction', instead, we
count harmonic primes up to some large bound.
The computation for a single prime in this case is fast and we go up to $50\cdot 10^5$,
extending Boyd's computation by a factor of $50$.
Our findings are summarized in the following theorem.

\begin{theorem}\label{thm:Jp}
$(i)$ For all primes $p\leq 16843$, the set $J_p$ is finite, with at most one exception, namely $p=1381$.

$(ii)$ There are $128594$ harmonic primes in the interval $[5,50\cdot 10^5]$,
corresponding to $\approx 36.89812\%$ of all primes in this range.

$(iii)$ There are no pairs $(p,n)$ with $p\leq 16843$, $p\neq 1381$,
for which $\nu_p(H_n)\geq 4$. If any such pair exists when $p=1381$, we must have $n\geq 1381^{3801}$.
\end{theorem}

The first point of Theorem \ref{thm:Jp} confirms Conjecture \ref{conj1}
for all primes $p\leq 16843$, with the exception of $1381$.
We did not complete the full enumeration of $J_{1831}$,
since we kept finding new elements all the way up to height $1381^{3800}$
(and in each of the last twenty $p$-adic intervals there are more than
$4.000$ elements, suggesting that we are far from completion).
A more precise version of point $(i)$ is stated in Observation \ref{obs1},
where we explain that $|J_p|\leq p^2$ for all the primes we examined
with four exceptions that satisfy instead the inequality $|J_p|\geq p^2(\log\log p)^2$.
This is in agreement with Boyd's quantitative version of Conjecture \ref{conj1}.
In Observation \ref{obs3} we also discuss the extinction time of $J_p$,
namely the largest power of $p$ needed to visit the whole set,
and compare it with the predictions from the model (see \cite[p.301]{Boyd}).

The second point in Theorem \ref{thm:Jp} (see Observation \ref{obs2})
is in agreement with the prediction
that harmonic primes should have density $e^{-1}=0.3678794411\dots$ among all primes
and hints at the correctness of Conjecture~\ref{conj2}.
Figure~\ref{Fig:harmonicprimes} shows the fluctuations around the value $e^{-1}$.

Finally, in the last point of Theorem \ref{thm:Jp}
we confirm that we never observe a $p$-adic valuation larger than $3$,
in agreement with Conjecture \ref{conj3}.
We found $21$ elements with valuation~$3$, see Observation \ref{obs4}
and Table \ref{Fig:valuation}.

Regarding progress towards a proof of either of Conjectures \ref{conj1}--\ref{conj3},
Sanna~\cite[Theorem 1.1]{S} proved that for any prime $p$ and any $x\geq 1$ we have
\[
|J_p \cap [1,x]| < 129 \, p^{\frac{2}{3}} \, x^{0.765}.
\]
Although not giving finiteness, this shows that $J_p$ has density zero
in the integers. Sanna's result has been improved by Wu and Chen~\cite[Theorem 1.1]{WC} to
\[
|J_p\cap [1,x]|\leq 3 \, x^{\frac{2}{3} + \frac{1}{25\log p}}.
\]
Bounds of this type have also been proved for harmonic numbers
of exponent greater than one by Altunta\c{s} \cite[Theorem A]{A}.
As for the possibility of having large $p$-adic valuation,
De Filpo and the first and third author showed that if
$p\nmid n$ and $\nu_p(H_n)$ equals $3$ or $4$ (resp.~\cite[Theorem 2.5 and 2.6]{CDG}),
then $\nu_p(H_{p^mn})$ grows linearly in $m$ before going down again to something $\leq 2$.
If we believe Conjecture \ref{conj3} is correct, then we should expect that the descent occurs immediately.
Our data confirms this, as we can see from Table~\ref{Fig:valuation}
where no two consecutive values of $m$ appear.

\noindent
All computations were made with \textsf{pari/gp} \cite{PARI};
source code is available online \cite{gc}.

\subsection*{Acknowledgments}
The second author is member of the INdAM group GNSAGA.
The third author is member of the INdAM group GNAMPA.

\section{Proof of Theorem \ref{thm:Jp}}

We wish to understand whether $J_p$ is finite or not,
what is the largest size $J_p$ can reach as $p$ varies, and what is
the largest $p$-adic valuation of its elements.
Our first step consists in splitting the integers in $p$-adic blocks
and checking the divisibility of harmonic numbers in each block:
let $m\geq 1$ and define the $m$-th $p$-adic block of $J_p$ as
\[
J_{p,m} := J_p \cap [p^{m-1},p^{m}-1].
\]
Clearly, $J_p$ is the union of the sets $J_{p,m}$ as $m$ varies.
Moreover, as explained by Boyd in \cite[\S 3]{Boyd},
$J_p$ has a recursive structure, so that its $m$-th block
can be obtained from the previous one,
provided we understand the latter sufficiently well.
To see this, let $k\in[0,p-1]$ and set $H_0=0$.
By \cite[Lemma 3.1]{Boyd} we have
\begin{equation}\label{1411:eq001}
H_{pn+k} = H_{pn} + H_{k} + O(p) = \frac{H_n}{p} + H_k + O(p).
\end{equation}
Here and in the rest of the paper we use the convention that
something is $O(p^s)$ if it is divisible by $p^s$.
Therefore, if we know $H_n$ up to an error $O(p^2)$ for all $n\in[p^{m-1},p^m-1]$,
as well as the value of $H_k$ up to an error $O(p)$ for all $k\in[0,p-1]$,
we can determine if $H_{pn+k}$ is a $p$-adic integer
and if $\nu_p(H_{pn+k})\geq 1$, for all integers $pn+k\in[p^{m},p^{m+1}-1]$.
In particular, \eqref{1411:eq001} implies that if $\nu_p(H_n)\leq 0$
for all integers $n$ in a given $p$-adic block,
then for all integers in the next block we have again $\nu_p(H_{n})\leq 0$.
In turn, this shows that the finiteness of $J_p$ is equivalent with showing
that $J_{p,m}=\emptyset$ for some $m\geq 1$, i.e.~that eventually there is an empty block.

By the above discussion it follows that the elements of $J_p$ can be arranged in a tree,
where the nodes at level $m$ are those $n$ in the interval $[p^{m-1},p^m-1]$
with $\nu_p(H_n)\geq 1$ and for every integer $k\in[0,p-1]$
there is an edge from $n$ to $pn+k$ if and only if $H_n\equiv -p H_k\bmod{p^2}$.
The set of residues $R=\{H_k\bmod{p}\}$ 
plays an important role in the structure of such a tree.
Assuming that the elements in $R$ are essentially randomly distributed,
Boyd \cite[\S 3 and \S 6--\S 7]{Boyd} gave a heuristic argument that suggests
that for every fixed $\epsilon>0$ we should have
\[
|J_p| \ll_\epsilon p^2(\log\log p)^{2+\epsilon}
\]
for all primes, although there should be infinitely
many primes for which
\begin{equation}\label{Jp:lowerbound}
|J_p| > p^2(\log\log p)^2.
\end{equation}
Boyd computed $J_p$ for all primes $p<550$ and his results were consistent with these
predictions. In particular, he found that $|J_{11}|=638>11^2$, giving one instance of \eqref{Jp:lowerbound}.
When $p=83,127$ and $397$, he could not determine the set $J_p$ in full
but obtained lower bounds on $|J_p|$ by looking at $p$-adic blocks $J_{p,m}$ with $m\leq 100$,
see \cite[p.288]{Boyd}.
We complete the computation for these three primes and go further, up to the first Wolstenholme prime,
obtaining the following.

\begin{obs}\label{obs1}
For all $5\leq p\leq 16843$ we have $|J_p|\leq p^2$, unless
$p=11,83,397$ or $1381$, in which case we have
\begin{equation}\label{Jp:large}
|J_{11}|=638,\quad
|J_{83}|=43038,\quad
|J_{397}|=701533,\quad
|J_{1381}|\geq 7521563.
\end{equation}
\end{obs}
The prime $p=127$ completes with precision $m=146$ and gives $|J_{127}|=3515$.
When $p=1381$ we could not complete the determination of $J_p$ (we reached precision $3800$)
and that is why we only have a lower bound in \eqref{Jp:large}.

\begin{figure}[!ht]
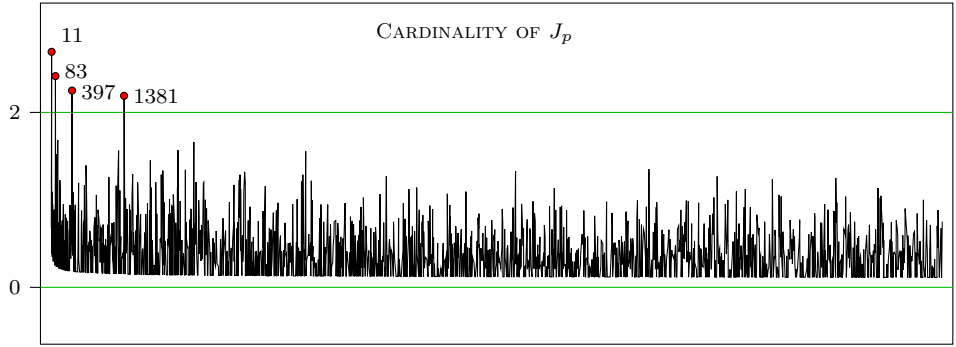

\begin{center}

\end{center}
\captionsetup{width=0.85\textwidth}
\caption{Cardinality of $J_p$ on a logarithmic scale.
On the horizontal axis we have $5\leq p\leq 16843$
and on the vertical axis the quantity $\log |J_p|/\log p$.
The profile on the bottom corresponds to the curve $\log 3/\log p$
associated with harmonic primes for which $|J_p|=3$.}
\end{figure}

One can also look at how the cardinalities $|J_p|$ distribute
as $p$ varies. They certainly do not distribute uniformly
but rather tend to favour small numbers. For instance, more than
sixty percent of all primes $p\leq 16843$ have $|J_p|\leq 31$
and about $36$ percent have $|J_p|=3$.
In Table \ref{Fig:Jpbins} the distribution among the
observed cardinalities up to $31$ is given. Curiously,
in this range not all integers are observed.
For instance, there is no prime $p\leq 16843$ with $|J_p|=5$.
Another visible feature is that most observed integers are odd.
This can partly be explained by Boyd's probabilistic model:
apart from the set $J_{p,1}$ which often contains the single element $p-1$,
the model predicts that at each successive level $J_{p,m}$ an even number
of elements is generated \cite[\S 6.2]{Boyd}, making the total count odd.
We indeed observe such a parity phenomenon in most levels.
Nevertheless, there are some primes with $|J_p|$ even, too.

\begin{table}[!ht]
\begin{center}
\scalebox{0.85}{%
\renewcommand\arraystretch{1.9}
\begin{tabular}{*{14}{|c}|}
\hline
\hline
\multicolumn{14}{c}{\textsc{Distribution of $|J_p|$}}\\
\hline
\hline
\rowcolor[gray]{0.85}3&7&9&11&13&15&17&19&21&23&25&27&29&31\\
\hline
\rowcolor[rgb]{.7,.95,.9}706&99&57&48&44&36&38&35&31&39&25&24&26&33\\
\hline\hline
\rowcolor[rgb]{.95,0.8,.9}36.35&5.09&2.93&2.47&2.26&1.85&1.95&1.80&1.59&2.00&1.28&1.23&1.33&1.69\\
\hline
\hline
\end{tabular}
}
\normalsize
\end{center}
\captionsetup{width=0.9\textwidth}
\caption{For $3\leq N\leq 31$, count of primes $5\leq p\leq 16843$
with $|J_p|=N$ and corresponding percentage of the total (the last digit is rounded down).
The values $18,20,24,26$ occur exactly once and are omitted. No other $N\leq 30$ appears.
Values above $31$ appear less than $19$ times each (less than $1\%$ of the total)
and are omitted.}\label{Fig:Jpbins}
\end{table}

The case $|J_p|=3$ is special, since for every $p\geq 5$ Eswarathasan and Levine \cite{EL}
showed that $J_p$ contains $p-1,p^2-1,p^2-p$ and therefore $|J_p|\geq 3$.
As we explained in the introduction, they called `harmonic primes'
those primes for which equality holds.
Table \ref{Fig:Jpbins} shows that out of $1942$
primes in the interval $[5,16843]$, about $36.35$ percent are harmonic.
Boyd's model predicts that harmonic primes should have density $e^{-1}=0.36787944\dots$
among all primes \cite[\S 4]{Boyd}. He computed harmonic primes up to $10^{5}$,
which agreed with such a prediction, although he writes that `the number of harmonic primes
in a given interval is perhaps somewhat higher than expected'.

\begin{figure}[!ht]
\begin{center}
\begin{tikzpicture}[scale=1.0,x=6pt,y=500pt]
\def\harmonicdatalist{
0/0.36430,
1/0.36205,
2/0.39674,
3/0.37161,
4/0.37742,
5/0.37338,
6/0.40319,
7/0.37916,
8/0.38813,
9/0.36746,
10/0.36469,
11/0.39151,
12/0.36364,
13/0.35840,
14/0.37709,
15/0.36527,
16/0.36978,
17/0.36568,
18/0.38768,
19/0.39312,
20/0.34508,
21/0.38964,
22/0.40659,
23/0.34184,
24/0.32685,
25/0.35183,
26/0.38385,
27/0.38734,
28/0.36742,
29/0.37904,
30/0.37484,
31/0.40842,
32/0.35678,
33/0.33933,
34/0.37484,
35/0.35256,
36/0.32810,
37/0.35476,
38/0.38462,
39/0.34174,
40/0.36605,
41/0.38918,
42/0.37047,
43/0.36329,
44/0.39608,
45/0.36835,
46/0.36601,
47/0.35038,
48/0.36794,
49/0.36788};

    \draw
	   (-0.25,0.367879) node[left] {$\frac{1}{e}$}
	   (0,0.32685)--(-0.4,0.32685) node[left] {\footnotesize$0.3268$}
	   (0,0.40842)--(-0.4,0.40842) node[left] {\footnotesize$0.4084$}
	   (50.4,0.33) node[right] {\phantom{0.0000}};

    \draw (25,0.467) node {\footnotesize\textsc{Count of harmonic primes}};

	\begin{scope}
	
    \path[clip] (0,0.3) rectangle (50,0.45);
	\draw (0,0) rectangle (50,0.78);
	\draw[green,thick] (0,0.367879)--(50,0.367879);
	
    \foreach \x/\y in \harmonicdatalist
		\draw[fill=gray!20,opacity=0.35] (\x,0) rectangle ({\x+1},\y);

    \draw[fill] (0,0.36430)
		\foreach \x/\y in \harmonicdatalist
			{-- ({\x+0.5},\y) circle [radius=0.75pt] };

    \end{scope}

    \draw (0,0.3) rectangle (50,0.45);

\end{tikzpicture}
\end{center}

\begin{center}
\begin{tikzpicture}[scale=1.0,x=6pt,y=1500pt]
\def\harmonicdatalist{
0/0.3777,
1/0.3736,
2/0.3679,
3/0.3618,
4/0.3705,
5/0.3755,
6/0.3590,
7/0.3629,
8/0.3796,
9/0.3576,
10/0.3713,
11/0.3634,
12/0.3710,
13/0.3744,
14/0.3611,
15/0.3592,
16/0.3651,
17/0.3610,
18/0.3719,
19/0.3709,
20/0.3658,
21/0.3744,
22/0.3730,
23/0.3692,
24/0.3679,
25/0.3659,
26/0.3622,
27/0.3734,
28/0.3664,
29/0.3705,
30/0.3614,
31/0.3697,
32/0.3738,
33/0.3629,
34/0.3762,
35/0.3659,
36/0.3652,
37/0.3701,
38/0.3579,
39/0.3737,
40/0.3671,
41/0.3722,
42/0.3834,
43/0.3734,
44/0.3693,
45/0.3761,
46/0.3816,
47/0.3680,
48/0.3603,
49/0.3762};

    \draw
	   (-0.25,0.367879) node[left] {$\frac{1}{e}$}
	   (0,0.3576)--(-0.4,0.3576) node[left] {\footnotesize$0.3576$}
	   (0,0.3834)--(-0.4,0.3834) node[left] {\footnotesize$0.3834$}
	   (50.4,0.33) node[right] {\phantom{0.0000}};

	\begin{scope}

    \path[clip,draw] (0,0.345) rectangle (50,0.395);
	\draw[green,thick] (0,0.367879)--(50,0.367879);

    \foreach \x/\y in \harmonicdatalist
		\draw[fill=gray!20,opacity=0.35] (\x,0) rectangle ({\x+1},\y);

    \draw[fill] (0,0.36430)
		\foreach \x/\y in \harmonicdatalist
			{-- ({\x+0.5},\y) circle [radius=0.75pt] };
	
    \end{scope}

\end{tikzpicture}
\end{center}
\captionsetup{width=0.865\textwidth}
\caption{Count of harmonic primes in $50$ intervals of size $10^4$ (top) and of size $10^5$ (bottom).
In the top part, the first 10 columns correspond to \cite[Table~1]{Boyd}.
The value $e^{-1}\approx0.367879$ is the density predicted by Boyd's probabilistic model.
}\label{Fig:harmonicprimes}
\end{figure}
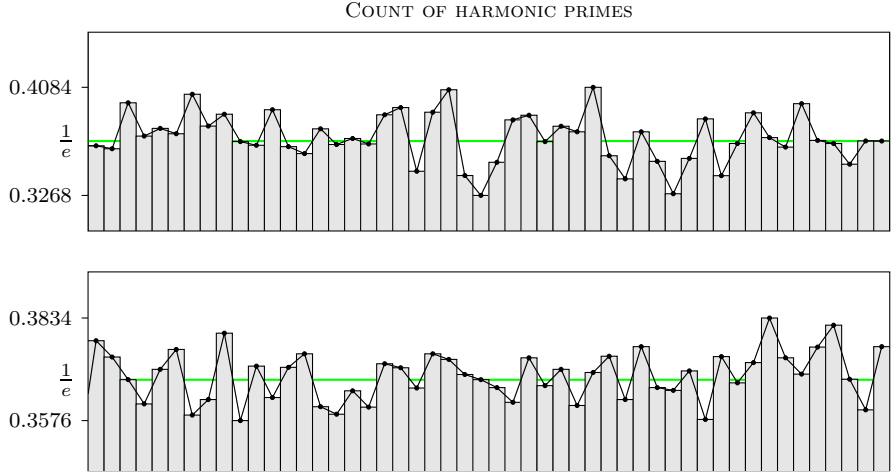

We extend Boyd's computation to primes up to $50\cdot 10^5$
and in Figure \ref{Fig:harmonicprimes} we plot the ratio
of harmonic primes in fifty intervals of size $10^4$ (top) and of size $10^5$ (bottom)
over all primes in the same interval.
There are fluctuations around the value $e^{-1}$ but the numerics
seem very convincing that this should be the correct density.
For instance, in the interval $[490000,500000]$
the fit is so accurate that we find
$284$ harmonic primes out of $772$ primes, for a ratio of
\[
\frac{284}{772} = 0.367875647\dots
\]
which agrees with $e^{-1}$ to the fifth decimal digit.
As for the total count, we have the following.

\begin{obs}\label{obs2}
Out of $348511$ primes $p\in[5,50\cdot 10^5]$, $128594$ are harmonic, giving a ratio $128594/348511\approx 0.3689812$.
\end{obs}

Returning to \eqref{1411:eq001},
let us finish to explain how the set $J_p$ is computed.
From \eqref{1411:eq001}, we see that $H_{pn}-H_n/p$ is a $p$-adic integer
with valuation at least one. More is true: there exists a sequence of
$p$-adic numbers $\{c_k\}_{k\geq 1}$ such that for all integers $n\geq 1$ we have
\begin{equation}\label{1209:eq001}
H_{pn}-\frac{1}{p}H_n = \sum_{k=1}^\infty c_k p^{2k} n^{2k}.
\end{equation}
This is proved in \cite[Theorem 5.2]{Boyd}.
We use \eqref{1209:eq001} to compute harmonic numbers as follows.
First, we compute the numbers $b_n=H_{pn}-H_n/p$
for $n=1,\dots,N$, to a $p$-adic precision $s$, and then solve the linear system
\begin{equation}\label{1209:eq002}
\sum_{k=1}^{N} c_k p^{2k} n^{2k} = b_n +O(p^s)
\end{equation}
in $c_1p^2,\dots,c_Np^{2N}$. The matrix of this system is the Vandermonde matrix $V$,
whose inverse satisfies $\nu_p(V^{-1})>-2N/(p-1)$, so that the unknown $c_kp^{2k}$
are obtained to precision $s-2N/(p-1)$. For fixed $N$,
the sum on the left in \eqref{1209:eq002} represents $b_n$ to precision
\cite[Remark 3]{Boyd}
\[
s = \min_{k>N} (\nu_p(c_k)+2k) \geq 2N+2-[\log_p(N+1)].
\]

An algorithm to calculate $J_p$ starts by computing $b_n=H_{pn}-H_n/p$
for $n\leq N$ and finding the coefficients $c_k'=c_kp^{2k}$
to precision $s\geq 2N+2-[\log_p(N+1)]$ as explained above.
In the process, one will have computed $H_n$ for $1\leq n\leq p-1$
to precision at least $s$ and hence will know $J_{p,1}$.
Then, once we have the elements at a given level $J_{p,m}$
to a precision $r\leq s$, one computes $H_{pn}$
from \eqref{1209:eq001} to precision $r-1$ and then computes
$H_{pn+k}=H_{pn+k-1}+1/(pn+k)$ for $k=1,\dots,p-1$,
thus obtaining $J_{p,m+1}$ to precision $r-1$.

Notice that the precision decreases by one at each new level and so we can
calculate elements in $J_p$ up to the $s$-th block $J_{p,s}$. If this set is empty,
then we are done and we have found all elements in $J_p$, which is finite.
If $J_{p,s}$ is not empty, we begin the computation again
with a larger value of $N$.

As pointed out in \cite[Section 5]{Boyd},
this method is faster than the `naive' method that goes by
computing harmonic numbers with the recursion $H_n=H_{n-1}+1/n$.
For instance, the naive method could not complete  in \cite{Boyd}
the full determination of $J_{11}$, which has 638 elements
and contains integers as large as $11^{30}$, whereas the $p$-adic
method described above succeeds and can go much further than that.

When running the algorithm with precision $s$, if $J_{p,s}\neq \emptyset$
we need to go back to the beginning and repeat the computation with a higher
precision. To speed up successive computations, we observe that at each level
not all nodes have children and so it is not necessary to compute
all elements in $J_p$, but only those that have descendants in the $s$-th block $J_{p,s}$.
This produces quite a bit of time and memory saving. For instance, when $p=1381$,
we find that $J_{1381,s}$ is non-empty for all $s\leq 3800$. If we look at intermediate levels,
we notice that $|J_{1381,1663}|=2501$ but only one element in this block has
descendants all the way down to $J_{1381,3800}$. Therefore, when running the algorithm
for any $s\geq 3800$, it suffices to calculate one element in each $J_{1381,r}$ for all $r\leq 1663$,
for a total of $1663$ harmonic numbers instead of $|J_{1381,1}\cup\cdots\cup J_{1381,1663}|=1860315$
elements.

\begin{figure}[!hb]
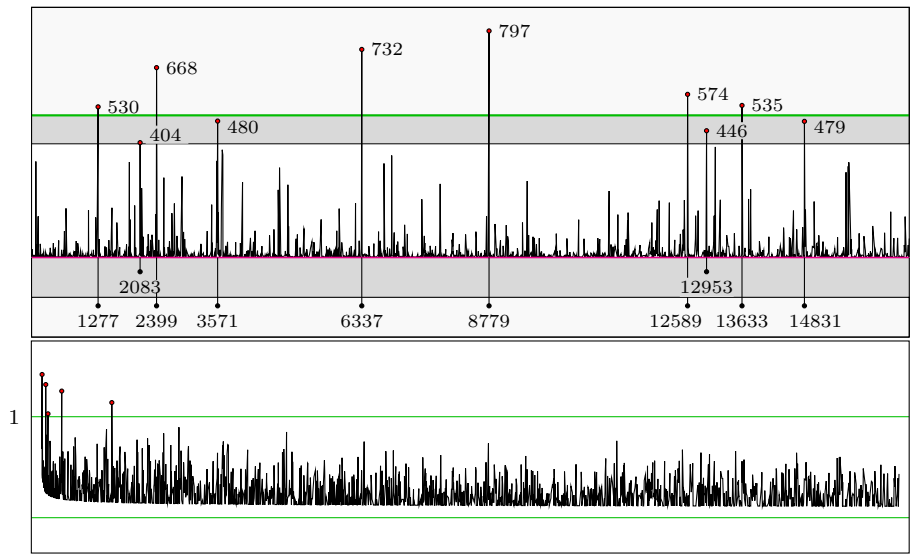

\begin{center}

\end{center}
\captionsetup{width=0.81\textwidth}
\caption{For primes $5\leq p\leq 16843$, we plot the extinction time $M_p$ (top figure; $p=397,1381,2699,4813,11299$ are omitted)
and in logarithmic scale we plot $\log M_p/\log p$ (bottom figure, including all primes).}\label{Fig:precision}
\end{figure}

In Figure \ref{Fig:precision} we plot the precision required to
compute $J_p$ for all primes $p\leq 16843$. This is sometimes referred to
as the `extinction time' for the branching process associated to $J_p$.
Similarly as with the cardinality, the random model predicts that
the extinction time should always be $O(p(\log\log p)^{1+\epsilon})$,
but there should be infinitely many primes with extinction time
larger than $p\log\log p$.

We see from Figure \ref{Fig:precision} that the extinction time is
indeed often large, say larger than 400. However,
very few primes have an extinction time as large as $p\log\log p$
and it does not come as a surprise that they are essentially the same ones
for which the cardinality is exceptionally large.
In fact, in the top part of the figure the primes $p=397$ and $1381$
are omitted, since their extinction time is way higher than all other primes
(respectively $1814$ and more than $3801$).
We also omit the primes $2699$, $4813$ and $11299$, whose extinction times are
respectively $1186$, $1336$ and $1214$.
The peaks in the bottom part of Figure \ref{Fig:precision}
correspond to $p=11,83,127,397$ and $1381$.
Summarizing, we have the following.

\begin{obs}\label{obs3}
The extinction time $M_p$ satisfies $M_p\leq p$ for all primes $p\leq 16843$,
with the following exceptions:
\[
M_{11}=30,\quad
M_{83}=339,\quad
M_{127}=146,\quad
M_{397}=1815,\quad
M_{1381}\geq 3801.
\]
\end{obs}

Finally, we conclude by discussing harmonic numbers with large $p$-adic valuation.
For an integer $n$ to be in $J_p$ we must have $\nu_p(H_n)\geq 1$
and computations reveal many integers for which $\nu_p(H_n)=2$.
Based on his model, Boyd conjectured that there are primes $p$ for which the
number of $n$ such that $\nu_p(H_n)=3$ is arbitrarily large,
but of order between $(\log\log p)^2$ and $(\log\log p)^{2+\epsilon}$.
On the other hand, he conjectured that $\nu_p(H_n)\geq 4$ never occurs.
In his work, he found no element with valuation $4$ or higher,
and only five instances of valuation 3: four when $p=11$
and one when $p=83$.
With the new data at hand, we have the following.

\begin{obs}\label{obs4}
$(i)$ For a given prime $5\leq p\leq 16843$,
the number of integers $n$ such that $\nu_p(H_n)=2$
can be as large as $5314$. More precisely,
\[
\max_{5\leq p\leq 16843} |\{n\in J_p:\nu_p(H_n)=2\}|
=
|\{n\in J_{1381}:\nu_{1381}(H_n)=2\}| \geq 5314.
\]
The second largest value is $1760$, which is attained when $p=397$.

\noindent$(ii)$
There are 21 pairs $(p,n)$ with $5\leq p\leq 16843$ for which $\nu_p(H_n)=3$
and the integers~$n$ appear in the
$p$-adic intervals $J_{p,m}$ described in Table \ref{Fig:valuation}.
When $p=1381$, possible additional occurrences
must have $n\geq 1381^{3801}$.

\noindent$(iii)$
There are no pairs $(p,n)$ with $5\leq p\leq 16843$, $p\neq 1381$,
for which $\nu_p(H_n)\geq 4$. If any such pair exists when $p=1381$, we must have $n\geq 1381^{3801}$.
\end{obs}

\begin{table}[!ht]
\begin{center}
\renewcommand\arraystretch{1.2}
\begin{tabular}{c|c|l|c}
\hline
& $p$ & $m$ &\\
\hline
&$11$ & $3,4,4,18$ &\\
\hline
&$83$ & $63,108,108,131,161,207,213,243,246,291,294$ &\\
\hline
&$397$ & $567$ &\\
\hline
&$1381$ & $1519,2572,2951,3211,3726$ &\\
\hline
\end{tabular}
\end{center}
\captionsetup{width=0.9\textwidth}
\caption{For each prime $p=11,83,397,1381$, elements with valuation three
are found in the intervals $J_{p,m}$ with $m$ as listed on the right.}\label{Fig:valuation}
\end{table}

Notice that, when $p=83$, the new integers we found
with valuation three are all larger than $p^{107}$.
Since Boyd computed $J_p$ up to $p^{100}$, this explains why they do not appear
in his work.


\end{document}